

\input epsf.tex

\def\2{{1\over 2}}

\def\d{\delta}
\def\a{\alpha}
\def\b{\beta}
\def\g{\gamma}

\def\s{\sigma}
\def\e{\epsilon}
\def\l{\lambda}

\def\fun#1#2#3{#1\colon #2\rightarrow #3}

\def\frac#1#2{{{#1} \over {#2}}}

\def\st{\;\colon\;}
\def\tends{\rightarrow}

\def\dr{ {\rm d} }

\def\R{{\bf R}}
\def\N{{\bf N}}

\def\thm#1{\vskip 1 pc\noindent{\bf Theorem #1.\quad}\sl}
\def\lem#1{\vskip 1 pc\noindent{\bf Lemma #1.\quad}\sl}
\def\prop#1{\vskip 1 pc\noindent{\bf Proposition #1.\quad}\sl}

\def\proof{\rm\vskip 1 pc\noindent{\bf Proof.\quad}}
\def\fin{\par\hfill $\backslash\backslash\backslash$\vskip 1 pc}
\def\txt#1{\quad\hbox{#1}\quad}

\def\L{{\cal L}}

\def\s{\sigma}

\def\r{\rho}

\def\2{\frac{1}{2}}
\def\inn#1#2{{\langle #1 ,#2\rangle}}

\def\part{{\partial_{x}}}

\def\pprime{{{}^\prime{}^\prime}}
\def\diam{{\rm diam}}
\def\tr{{{}^{t}}}

\def\ec{{\cal E}}

\def\dc{{\cal D}}
\def\oc{{\cal O}}
\def\pc{{\cal P}}

\def\mc{{\cal M}}



\baselineskip= 17.2pt plus 0.6pt
\font\titlefont=cmr17
\centerline{\titlefont Cheeger's energy}
\vskip 1 pc
\centerline{\titlefont on the harmonic Sierpinski gasket}
\vskip 4pc
\font\titlefont=cmr12
\centerline{         \titlefont {Ugo Bessi}\footnote*{{\rm 
Dipartimento di Matematica, Universit\`a\ Roma Tre, Largo S. 
Leonardo Murialdo, 00146 Roma, Italy.}}   }{}\footnote{}{
{{\tt email:} {\tt bessi@matrm3.mat.uniroma3.it} Work partially supported by the PRIN2009 grant "Critical Point Theory and Perturbative Methods for Nonlinear Differential Equations}} 
\vskip 0.5 pc
 
\par
\vskip 2pc
\centerline{\bf Abstract}

Koskela and Zhou have proven that, on the harmonic Sierpinski gasket with Kusuoka's measure, the "natural" Dirichlet form coincides with Cheeger's energy. We give a different proof of this result, which uses the properties of the Lyapounov exponent of the gasket.

\vskip 2 pc
\centerline{\bf  Introduction}
\vskip 1 pc

This paper deals with the harmonic Sierpinski gasket $S$ endowed with Kusuoka's measure $\kappa$; we refer the reader to sections 1 and 2 below for the definitions and properties of these objects. 

Let $P_v$ denote the orthogonal projection in $\R^2$ along the vector $v$: 
$$P_v w=\frac{(w,v)}{||v||^2}\cdot v   .  $$
As we recall in section 2 below, there is a Borel vector field $\fun{v}{S}{\R^2}$ such that, if one defines 
$$\fun{\ec}{C^1(\R^2,\R)\times C^1(\R^2,\R)}{\R}$$
$$\ec(\phi,\psi)=\int_S(\nabla\phi,P_{v(x)}\nabla\psi)\dr\kappa(x)   \eqno (1)  $$
then $\ec$ extends to a local, self-similar Dirichlet form densely defined on 
$L^2(S,\kappa)$. Heuristically, this means that $v_x$ is the "tangent space at $x$" seen by the Brownian motion on $S$. As a remark on notation, we shall denote the vector field at $x\in S$ both by $v(x)$ and by $v_x$. 

There is another quadratic form densely defined on $L^2(S,\kappa)$, and it is Cheeger's energy ([2], [4]); in [12] (see also [8]), the following theorem is proven.

\thm{1} The Dirichlet form $\ec$ coincides with the double of Cheeger's energy; in symbols,    
$$Ch(u)=\2\ec(u,u)    \eqno (2)$$
for all $u\in\dc(Ch)=\dc(\ec)$. 

\vskip 1pc

\rm

Morally, this means that $v_x$ is also the "tangent space at $x$" seen by large families of absolutely continuous curves, the so-called test plans with bounded deformation. 

The proof of [12] is rather delicate; the aim of this paper is to provide a different proof which links Cheeger's energy to the Lyapounov exponent of the dynamical system underlying 
$S$.  Indeed, a third way of looking at the "tangent space" of $S$ is the following: there is an expansive map $\fun{F}{\R^2}{\R^2}$ which sends $S$ into itself and whose derivative is defined $\kappa$-a. e. on $S$. Heuristically, the vector field $v(x)$ of (1) coincides with the "less expansive" direction of $F$ and it determines the only viable velocity for curves in $S$; more precisely, we shall prove that for "almost" all absolutely continuous curves $\fun{\g}{[0,1]}{S}$ which are differentiable at $t\in(0,1)$, we have (supposing, which is always possible, that 
$||v(x)||\equiv 1$)
$$\dot\g_t=\pm ||\dot\g_t||\cdot v(\g_t)  .  $$
If $\phi\in C^1(\R^2,\R)$, the formula above implies, practically by definition of Cheeger's derivative $|D\phi|_w$, that 
$$|D\phi|_w(x)\le |(\nabla \phi,v(x))|
\txt{for $\kappa$-a. e. $x\in S$}  $$
which is half of the proof of formula (2).  

For the opposite inequality, we  use two facts. First, we prove that in Cheeger's definition of the energy (formula (1.13) below) we can restrict to $C^1$ functions. Second, using again the Lyapounov exponent we show that, if $f\in C^1$, then its local Lipschitz constant on $S$, usually denoted by $Lip_a(f,x)$, is larger than $|(\nabla f(x),v(x))|$.

The paper is organised as follows. In section 1 below we recall the main definitions and properties of the Harmonic Sierpinski gasket and Cheeger's energy. In section 2 we introduce Kusuoka's measure, following the approach of [3]. In section 3, we study the shape of smaller and smaller cells of the gasket: heuristically, they are "skinny" triangles with the long side aligned with $v(x)$. We shall also see how Kusuoka's measure distributes inside these triangles. Theorem 1 is proven in section 4. 

\vskip 1pc

\noindent{\bf Acknowledgements.} The author would like to thank the referee for the useful comments.

\vskip 2pc

\centerline{\bf\S 1}

\centerline{\bf The gasket and Cheeger's energy}

\vskip 1pc

\noindent{\bf The gasket.} We recall the definition of the harmonic Sierpinski gasket; we follow [8]. We consider the three linear contractions of $\R^2$
$$T_1=
\left(
\matrix{
\frac{3}{5} &0\cr
0&\frac{1}{5}
}
\right),\qquad 
T_2=
\left(
\matrix{
\frac{3}{10} &\frac{\sqrt 3}{10}\cr
\frac{\sqrt 3}{10}&\frac{1}{2}
}
\right),\qquad 
T_3=
\left(
\matrix{
\frac{3}{10} &-\frac{\sqrt 3}{10}\cr
-\frac{\sqrt 3}{10}&\frac{1}{2}
}
\right)    $$
and the three vertices of an equilateral triangle $\bar A\bar B\bar C$: 
$$\bar A=(0,0),\qquad
\bar B=\left(  1,\frac{1}{\sqrt 3}\right  ),\qquad
\bar C=\left(  1,-\frac{1}{\sqrt 3}  \right)   .   $$
We define the three affine contractions of $\R^2$ 
$$\psi_1(x)=\bar A+T_1(x-\bar A),\qquad
\psi_2(x)=\bar B+T_2(x-\bar B),\qquad
\psi_3(x)=\bar C+T_3(x-\bar C)   .   \eqno (1.1)$$
It turns out (theorem 1.1.7 of [9]) that there is a unique compact set $S\subset\R^2$, which is called the harmonic Sierpinski gasket, such that 
$$S=\bigcup_{i=1}^3\psi_i(S)  .  \eqno (1.2)$$
The set $S$ is contained in the equilateral triangle  $\bar A\bar B\bar C$; since the three vertices $\bar A$, $\bar B$ and $\bar C$ are fixed points of $\psi_1$, $\psi_2$ and 
$\psi_3$ respectively, they are contained in $S$; they are often called the boundary of $S$; we refer the reader to [15] for this terminology and other interesting facts on the gasket. 

The three maps $\psi_1$, $\psi_2$ and $\psi_3$ are the three branches of the inverse of an expansive map $\fun{F}{S}{S}$. More precisely, let us consider the three triangles 
$\psi_i(\bar A\bar B\bar C)$ for $i=1,2,3$; they are depicted on the right of figure 1 below, where $\psi_1(\bar A\bar B\bar C)=\bar Abc$, $\psi_2(\bar A\bar B\bar C)=\bar Bca$ and 
$\psi_3(\bar A\bar B\bar C)=\bar Cba$.  Note that the sets $\psi_i(\bar A\bar B\bar C)$ intersect only at the vertices denoted by lowercase letters.

We can find three disjoint  open sets $\oc_1$, $\oc_2$ and $\oc_3$ such that $\oc_1$ contains the closed triangle $Abc$ save the two points $\{ b,c \}$, $\oc_2$ contains $Bca$ save $\{ c,a \}$ and $\oc_3$ contains $Cba$ save $\{ b,a \}$. We can also require that 
$\oc_j\subset\psi^{-1}(\oc_j)$.  We define
$$\fun{F}{
\bigcup_{i=1}^3 \oc_i
}{(\bar A\bar B\bar C)}  $$
by
$$F(x)=\psi_i^{-1}(x)  \txt{if}
x\in\oc_i  .  $$
This does not define $F$ on the three points $\{ a,b,c \}$; we set arbitrarily $F(a)=C$, 
$F(b)=A$ and $F(c)=B$. We can afford this arbitrariness since it is standard ([13], see [3] for an alternative proof) that $\{ a,b,c \}$ is a null set for Kusuoka's measure. 

A basic fact is that the dynamics of $F$ on $S$ can be coded. Namely, let us consider the space of sequences 
$$\Sigma\colon=\left\{
\{ x_i \}_{i\ge 0}\st x_i\in(1,2,3)\qquad\forall i\ge 0
\right\}  $$
with the product topology. Let $\eta\in(0,1)$ be the common Lipschitz constant of $\psi_1$, $\psi_2$ and $\psi_3$ (by (1.1) we gather that $\eta=\frac{3}{5}$, but the precise value is immaterial). If $x=(x_0x_1\dots)\in\Sigma$, using the fact that 
$\diam(\bar A\bar B\bar C)=\frac{2}{\sqrt 3}$ we get that 
$$\diam(\psi_{x_0}\circ\dots\circ\psi_{x_l}(\bar A\bar B\bar C))\le
\frac{2}{\sqrt 3}\eta^{l+1}  .  \eqno (1.3)$$
Since $\bar A\bar B\bar C$ is a closed triangle, the sets 
$\psi_{x_0}\circ\dots\circ\psi_{x_l}(\bar A\bar B\bar C)$ are compact; formula (1.2) implies that $\psi_i(S)\subset S$ which in turn implies that for all $l\ge 1$ 
$$\psi_{x_0}\circ\dots\circ\psi_{x_l}\circ\psi_{x_{l+1}}(\bar A\bar B\bar C)
\subset
\psi_{x_0}\circ\dots\circ\psi_{x_l}(\bar A\bar B\bar C)  .  \eqno (1.4)$$
Together with (1.3) this implies that 
$$\bigcap_{l\ge 1}\psi_{x_0}\circ\dots\circ\psi_{x_l}(\bar A\bar B\bar C)$$
is a single point, which we call $\Phi(x)=\Phi(x_0x_1\dots)$. Formulas (1.3) and (1.4) imply easily that the map $\fun{\Phi}{\Sigma}{S}$ is continuous. 

Recall that the triangles $\psi_i(\bar A\bar B\bar C)$ intersect only at the vertices, which have the lowercase letters in the figure above; this easily implies that, if $x\in S$ is coded by $w_0w_1\dots$ and $w_0^\prime w_2^\prime\dots$ with $w_0\not=w_0^\prime$, then $x\in\{ a,b,c \}$. Iterating, we get that the points in $S$ which have multiple codings are of the type 
$$\psi_{x_0}\circ\dots\circ\psi_{x_l}(P)   
\txt{with}P\in(\bar A,\bar B,\bar C)  .  $$
Note that the points defined by the formula above are a countable set. It is easy to see that these points have at most two pre-images, i. e. that $\Phi$ is at most two to one. 

Next we note that, if $x\in S\setminus\{ B,C \}$, then $\psi_1(x)\in\oc_1$ and thus 
$F\circ\psi_1(x)=x$; the same argument holds for $\psi_2$  and $\psi_3$ yielding (1.5) below, while (1.6) is immediate.  
$$\left\{
\eqalign{
F\circ\psi_1(x)&=x\qquad\forall x\in G\setminus\{ B,C \}, \cr  
F\circ\psi_2(x)&=x\qquad\forall x\in G\setminus\{ A,C \}, \cr
F\circ\psi_3(x)&=x\qquad\forall x\in G\setminus\{ A,B \}. 
}
\right.
\eqno (1.5)$$
$$\psi_i\circ F(y)=y\qquad\forall y\in\oc_i  .  
\eqno (1.6)$$
By the chain rule, the last two formulas imply the two equalities below; since $D\psi_j$ is constant, we have omitted its argument.  
$$DF|_{\psi_j(x)}=(D\psi_j)^{-1}\qquad
\forall x\in G\setminus \{ A,B,C \},\qquad\forall j\in(1,2,3)   \eqno (1.7)$$
$$D\psi_j=(DF(y))^{-1}\qquad
\forall y\in \oc_j   .  \eqno (1.8)$$

Let us define the one-sided shift as 
$$\fun{\s}{\Sigma}{\Sigma}$$
$$\fun{\s}{\{ x_i \}_{i\ge 0}}{\{ x_{i+1} \}_{i\ge 0}}   .   $$
From (1.5) it follows easily that $\Phi$ conjugates $F$ and $\sigma$:
$$\Phi\circ\s(\{ x_i \}_{i\ge 0})=F\circ\Phi(\{ x_i \}_{i\ge 0})  \eqno (1.9)$$
save for the points $\Phi(x)\in S$ where $F\circ\Phi(x)$ is defined arbitrarily, i. e. 
$\{ a,b,c \}$. 


Let $x=(x_0x_1\dots)\in\Sigma$; we shall often use the notation  
$$\psi_{x,l}=\psi_{x_0\dots x_l}=
\psi_{x_0}\circ\dots\circ\psi_{x_l}   .  \eqno (1.10)$$
Having two ways to indicate the composition $\psi_{x_0}\circ\dots\circ\psi_{x_l}$ looks redundant, but we shall use both of them. 

We define the cell (or cylinder) $[x_0\dots x_l]$ as 
$$[x_0\dots x_l]=\psi_{x,l}(S)=
\psi_{x_0\dots x_l}(S)  .  \eqno (1.11)$$

\vskip 1pc

\noindent{\bf Absolutely continuous curves.} It is a standard fact ([8]) that $S$ is arcwise connected; this prompts us to recall a definition of [1].  

\vskip 1pc

\noindent{\bf Definition.} We say that $\g\in AC^2([0,1],S)$ if 

\noindent 1) $\g$ is absolutely continuous from $[0,1]$ to $\R^2$,

$$\int_0^1||\dot\g_s||^2\dr s<+\infty  \leqno 2)$$

\noindent 3)  $\g_s\in S$ for all $s\in[0,1]$. 

For $x,y\in S$ we define
$${\rm d}_{geod}^2(x,y)=\min
\left\{
\int_0^1||\dot\g_s||^2\dr s\st\g\in AC^2([0,1],S),\quad \g_0=x,\g_1=y
\right\}   .  $$
It turns out ([8], see also [19] and [10]) that the minimum is attained and that the minimal curves are of class $C^1$; they are geodesics of constant speed, i. e. for all 
$0\le s\le t\le 1$ we have 
$${\rm d}_{geod}(\g_s,\g_t)=(t-s){\rm d}_{geod}(x,y)  .  $$
By [8], the geodesic of constant speed connecting $x$ and $y$ is unique; we call it 
$\g_{x,y}$.  

Another fact (again [8]) is that the union 
$$M\colon =\g_{\bar A\bar B}([0,1])\cup\g_{\bar B\bar C}([0,1])\cup
\g_{\bar A\bar C}([0,1])  \eqno (1.12)  $$
is a closed set which divides the plane into a bounded and an unbounded component; $S$ is contained in the closure of the bounded component. 

\vskip 1pc

\noindent{\bf Test plans and Cheeger's energy.} A Borel probability measure on 
$C([0,1],S)$  is called a test plan if the two points below hold.

\noindent 1) $\pi$ concentrates on $AC^2([0,1],S)$.

\noindent 2) $\pi$ satisfies 
$$\int_{C([0,1],S)}\dr\pi(\g)\int_{0}^1||\dot\g_s||^2\dr s<+\infty  .  $$

Let us now fix a Borel probability measure $m$ on $S$; from section 3 onwards, $m$ will be Kusuoka's measure $\kappa$, but the definition we give below holds for a completely arbitrary $m$. For $t\in[0,1]$ let us denote by $e_t$ the evaluation map 
$$\fun{e_t}{C([0,1],S)}{S}$$
$$\fun{e_t}{\g}{\g_t}  .  $$
If $\pi$ is a test plan, let us define by $(e_t)_\sharp\pi$ the push-forward of $\pi$ by $e_t$; in other words, if $\fun{f}{S}{\R}$ is a bounded Borel function, then 
$$\int_S f(x)\dr[(e_t)_\sharp\pi](x)=
\int_{C([0,1],S)}f(\g_t)\dr\pi(\g)  .  $$
We say that the test plan $\pi$ has bounded deformation with respect to $m$ if for all 
$t\in[0,1]$ we have that 
$$(e_t)_\sharp\pi=\r_t m$$
and if there is $C>0$ such that 
$$||\r_t||_{L^\infty(m)}\le C\qquad
\forall t\in[0,1]  .  $$

We recall the definition of weak gradient from [2]. 





\vskip 1pc

\noindent{\bf Definition.} Let $f\in L^2(S,m)$; we say that $g\in L^2(S,m)$ is a weak gradient of $f$ if, for all test plans $\pi$ with bounded deformation with respect to $m$ and for 
$\pi$-a. e. curve 
$\g\in C([0,1],S)$ we have 
$$|f(\g_1)-f(\g_{0})|\le
\int_{0}^1 g(\g_s)||\dot\g_s||\dr s  .  $$
If $f$ admits a weak gradient $g\in L^2(S,m)$, then it admits a weak gradient $|Df|_w$ minimal in the following two senses: if $g\in L^2(S,m)$ is a weak gradient, then 
$$|||Df|_w||_{L^2(S,m)}\le ||g||_{L^2(S,m)}  $$ 
and
$$|Df|_w(x)\le g(x)\txt{$m$-a. e..}  $$

Cheeger's energy is  the function 
$$\fun{Ch}{L^2(S,m)}{[0,+\infty]}$$
defined by
$$Ch(f)=\2\int_S|Df|_w^2\dr m$$
if $f$ has a weak gradient in $L^2(S,m)$, and $Ch(f)=+\infty$ otherwise. 

Though there is a square in its definition, $Ch$ in general is not a quadratic form (see [2] for an example). However, if $m=\kappa$, then theorem 1 of the introduction implies that $Ch$ is quadratic.

We shall need an equivalent definition ([2], [4]) of $|Df|_w$. Namely, if 
$\fun{f}{S}{\R}$ is a function, we define 
$$|Df|(x)=\lim_{r\tends 0}
\sup_{y\in [B(x,r)\setminus\{ x \}]\cap S}
\frac{|f(y)-f(x)|}{||y-x||}    $$
where $||\cdot||$ is the euclidean norm on $\R^2$. 

It turns out that, if $f, |Df|_w\in L^2(S,m)$, then $|Df|_w$ is the $L^2$ function $g$ of minimal norm for which there is a sequence of Lipschitz functions $\fun{f_n}{S}{\R}$ such that 

\noindent 1) $f_n\tends f$ in $L^2(S,m)$ and

\noindent 2) $|Df_n|\tends g$ weakly in $L^2(S,m)$.

Points 1) and 2) above are Cheeger's characterisation of $|Df|_w$; this characterisation implies that  
$$Ch(f)=\inf\liminf_{n\tends+\infty}\2\int_S|Df_n|^2\dr m  $$
where the $\inf$ is over all sequences $\{ f_n \}_{n\ge 1}$ of Lipschitz functions such that 
$f_n\tends f$ in $L^2(S,m)$. 

In section 4 we shall use a different formula, which we found in chapter 3 of [5]. We define the local Lipschitz constant of $f$ by 
$$Lip_a(f,x)=
\lim_{r\tends 0}\sup\left\{
\frac{|f(y)-f(z)|}{||y-z||}\st z\not=y,\qquad z,y\in B(x,r)
\right\}  .  $$
Then,
$$Ch(f)=\inf\liminf_{n\tends+\infty}\2\int_S|Lip_a(f_n)|^2\dr m  \eqno (1.13)$$
where, again, the $\inf$ is over all sequences $\{ f_n \}_{n\ge 1}$ of Lipschitz functions which converge to $f$ in $L^2(S,m)$.

\vskip 2pc
\centerline{\bf \S 2}
\centerline{\bf Kusuoka's measure and fields of projections}

\vskip 1pc

In this section, we recall the definition and some properties of Kusuoka's measure ([13]); we follow the approach of [3]. 

Let $M^2$ denote the space of $2\times 2$ symmetric matrices; $M^2$ is a Hilbert space for the Hilbert-Schmidt inner product 
$$(A,B)_{HS}\colon=
tr(\tr A B)=tr(AB)$$
where $\tr A$ denotes the transpose of $A$; for the second equality, recall that $A$ is symmetric. 

We denote by $C(S,M^2)$ the space of continuous functions from $S$ to $M^2$; its dual is the space $\mc(S,M^2)$ of the $M^2$-valued Borel measures on $S$. 

We want to define the integral against $\mu\in\mc(S,M^2)$; in order to do this, we recall that the total variation of $\mu$ is a finite measure $||\mu||$ on the Borel sets of $S$. The polar decomposition of $\mu$ is given by 
$$\mu=M_x||\mu||$$
where $\fun{M}{S}{M^2}$ is a Borel field of symmetric matrices which satisfies 
$$||M_x||_{HS}=1
\txt{for $||\mu||$-a. e. $x\in S$.}  \eqno (2.1)$$
If $\fun{A}{S}{M^2}$ is Borel and $||A||_{HS}\in L^1(S,||\mu||)$, then by (2.1) and Cauchy-Schwarz we have that 

\noindent $(A_x,M_x)_{HS}\in L^1(S,||\mu||)$. Consequently, we can define the scalar 
$$\int_S(A_x,\dr\mu(x))_{HS}\colon =
\int_S(A_x,M_x)_{HS}\dr||\mu||(x)  .  $$
The duality coupling between $C(S,M^2)$ and $\mc(S,M^2)$ is given by 
$$\fun{\inn{\cdot}{\cdot}}{
C(S,M^2)\times\mc(S,M^2)
}{\R}  $$
$$\inn{A}{\mu}\colon =
\int_S(A_x,\dr\mu(x))_{HS}  .  $$
If $Q\in C(S,M^2)$ and $\mu\in\mc(S,M^2)$, we define the scalar measure 
$(Q,\mu)_{HS}$ by 
$$\int_S f(x)\dr(Q,\mu)_{HS}(x)\colon =
\int_S(fQ,\dr\mu)_{HS}    \eqno (2.2)$$
for all $f\in C(S,\R)$. In other words, $(Q,\mu)_{HS}=(Q_x,M_x)_{HS}\cdot ||\mu||$. 

If $\fun{v,w}{S}{\R^2}$ are Borel functions such that $||v||\cdot ||w||\in L^1(S,||\mu||)$ then, again by (2.1) and Cauchy-Schwarz, $(v_x,M_xw_x)\in L^1(S,||\mu||)$ and we can define 
$$\int_S(v_x,\dr\mu(x)w_x)\colon =
\int_S(v_x,M_xw_x)\dr||\mu||(x)  .  \eqno (2.3)$$

We say that $\mu\in\mc(S,M^2)$ is semipositive definite if $\mu(E)$ is a semipositive definite matrix for all Borel sets $E\subset S$. Lusin's theorem easily implies ([3]) that 
$\mu$ is semipositive definite if and only if 
$$\inn{A}{\mu}\ge 0   \eqno (2.4)$$
for all $A\in C(S,M^2)$ such that $A_x\ge 0$ for all $x\in S$. 

We denote by $\mc_+(S,M^2)$ the set of all semipositive definite measures of $\mc(S,M^2)$; by the characterisation of (2.4), $\mc_+(S,M^2)$ is a convex set of 
$\mc(S,M^2)$, closed for the  weak$\ast$ topology. 

Let now $Q\in C(S,M^2)$ be such that $Q_x$ is positive-definite for all $x\in S$; since $S$ is compact there is $\e>0$ such that 
$$Q_x\ge\e Id\qquad\forall x\in S  .  \eqno (2.5)$$
For such a $Q$ we define $\pc_Q(S,M^2)$ as the set of all $\mu\in\mc_+(S,M^2)$ such that 
$$\int_S(Q_x,\dr\mu(x))_{HS}=1  .  $$
As shown in [3], if $Q$ satisfies (2.5) there is $D_1(\e)>0$ such that, for all 
$\mu\in\mc_+(S,M^2)$, 
$$||\mu||\le D_1(\e)(Q,\mu)_{HS}   .  \eqno (2.6)$$
In other words, the Radon-Nykodim derivative of $||\mu||$ with respect to $(Q,\mu)_{HS}$ is bounded by $D_1(\e)$. As a consequence, if $\mu\in\pc_Q(S,M^2)$, then
$$||\mu||(S)\le D_1(\e)  .  $$
By the characterisation (2.4), $\pc_Q(S,M^2)$ is a convex subset of $\mc(S,M^2)$, closed for the weak$\ast$ topology; by the formula above, it is compact. 

Let $\psi_1, \psi_2, \psi_3$ be the affine maps of section 1; we define the Ruelle operator as 
$$\fun{\L}{C(S,M^2)}{C(S,M^2)}$$
$$(\L A)(x)\colon =\sum_{i=1}^3
\tr D\psi_i(x)\cdot A_{\psi_i(x)}\cdot  D\psi_i(x)  .  $$
The adjoint of $\L$ 
$$\fun{\L^\ast}{\mc(S,M^2)}{\mc(S,M^2)} $$
is defined by 
$$\inn{\L A}{\mu}=\inn{A}{\L^\ast\mu}  $$
for all $A\in C(S,M^2)$ and $\mu\in\mc(S,M^2)$. 

The following proposition is proven as in the scalar case, for which we refer to [16] and [20]; the details for the matrix case are in [3]. 

\prop{2.1} 1) The operator $\L$ has a simple, positive eigenvalue $\b>0$. Let 
$\tilde Q\in C(S,M^2)$ be an eigenfunction of $\b$; then, up to multiplying $\tilde Q$ by a scalar, we have 
$$\tilde Q_x=Id\qquad\forall x\in S. $$
In other words, the eigenspace of $\b$ is generated by a matrix field constantly equal to the identity. 

\noindent 2) There is a unique $\tau\in\pc_{Id}(S,M^2)$ such that 
$$\L^\ast\tau=\beta\tau  .  $$

\noindent 3) Let us define the scalar measure $\kappa\colon=(Id,\tau)_{HS}$ as in (2.2); then, $\kappa$ is a probability measure ergodic for the map $F$. Moreover, $\kappa$ is non-atomic. 

\noindent 4) For $f,g\in C^1(\R^2,\R)$ we define 
$$\ec(f,g)=\int_S(\nabla f(x),\dr\tau(x)\nabla g(x))    $$
where the integral has been defined in (2.3). Then, $\ec$ is self-similar; in other words, for the maps $\psi_i$ of section 1 and all $f,g\in C^1(\R^2,\R)$ we have that 
$$\ec(f,g)=\frac{1}{\b}\sum_{i=1}^3\ec(f\circ\psi_i,g\circ\psi_i)   .   $$

\noindent 5) The measure $\tau$ has the Gibbs property; in other words, with the notation of (1.10) and (1.11) we have that 
$$\tau[x_0\dots x_{l-1}]=
\frac{1}{\b^l}(D\psi_{x_0\dots x_{l-1}})\cdot
\tau(S)\cdot \tr(D\psi_{x_0\dots x_{l-1}})  .   $$
This formula further simplifies since (see for instance [3] or [8]) $\tau(S)=Id$. Since 
$\psi_{x_0\dots x_l}$ is affine, its derivative is constant; therefore, we haven't specified the point where we calculate it.

\rm

\vskip 1pc

Lemma 2.2 below gives an expression for the push-forward 
$(\psi_{j_0\dots j_{n-1}})_\sharp\tau$; the notation is that of (1.10).

\lem{2.2} Let $j_0,\dots, j_{n-1}\in(1,2,3)$; then,  
$$(\psi_{j_0\dots j_{n-1}})_\sharp\tau=
\b^n(D\psi_{j_0\dots j_{n-1}})^{-1}\cdot\tau|_{\psi_{j_0\dots j_{n-1}}(S)}\cdot\tr (D\psi_{j_0\dots j_{n-1}})^{-1}  .  
\eqno (2.7)$$
By $\tau|_{\psi_{j_0\dots j_{n-1}}(S)}$ we have denoted  the restriction of $\tau$ to 
$\psi_{j_0\dots j_{n-1}}(S)$; note that $(\psi_{j_0\dots j_{n-1}})_\sharp\tau$ is supported on $\psi_{j_0\dots j_{n-1}}(S)$ by definition. 

\proof Since the sets $[x_0\dots x_{l-1}]$ generate the Borel sets of $S$, it suffices to show that, for all $l\ge n\ge 0$ and $\{j_i\}_{i\ge 0},\{ x_i \}_{i\ge 0}\in\Sigma$, 
$$[(\psi_{j_0\dots j_{n-1}})_\sharp\tau][x_0\dots x_{l-1}]=
\b^n(D\psi_{j_0\dots j_{n-1}})^{-1}\cdot\tau ([x_0\dots x_{l-1}]\cap\psi_{j_0\dots j_{n-1}}(S) )\cdot\tr (D\psi_{j_0\dots j_{n-1}})^{-1}   .   
\eqno (2.8)$$
Since $(\psi_{j_0\dots j_{n-1}})_\sharp\tau$ is supported on the cell $[j_0\dots j_{n-1}]$, we can as well suppose that $x_0x_1\dots x_{n-1}=j_0j_1\dots j_{n-1}$: otherwise, we have zero on both sides of (2.8). Using this fact together with (1.10) and (1.11) we get the equality below. 
$$[x_{n}x_{n+1}\dots x_{l-1}]=
\psi_{j_0\dots j_{n-1}}^{-1}([x_0x_1\dots x_{l-1}])   .   \eqno (2.9)$$
By point 5) of proposition 2.1, we get the two equalities below:
$$\tau[\psi_{x_0\dots x_{l-1}}(S)]=
\frac{1}{\b^l}(D\psi_{x_0\dots x_{l-1}})\cdot\tau(S)\cdot\tr(D\psi_{x_0\dots x_{l-1}})  $$
$$\tau[\psi_{x_n\dots x_{l-1}}(S)]=
\frac{1}{\b^{l-n}}(D\psi_{x_n\dots x_{l-1}})\cdot\tau(S)\cdot\tr(D\psi_{x_n\dots x_{l-1}}) . $$
By the chain rule, the last two formulas imply that 
$$\tau[\psi_{x_n\dots x_{l-1}}(S)]=
\b^n(D\psi_{x_0\dots x_{n-1}})^{-1}\cdot
\tau[\psi_{x_0\dots x_{l-1}}(S)]\cdot
\tr(D\psi_{x_0\dots x_{n-1}})^{-1}   .  \eqno (2.10)$$
The first equality below is the definition of push-forward; the second one comes from (2.9), the third one from (2.10) and the fact that $x_0x_1\dots x_{n-1}=j_0j_1\dots j_{n-1}$. 
$$[(\psi_{j_0\dots j_{n-1}})_\sharp\tau][x_0\dots x_{l-1}]=
\tau(\psi_{j_0\dots j_{n-1}}^{-1}([x_0\dots x_{l-1}]))=
\tau[x_{n}\dots x_{l-1}]=$$
$$\b^n (D\psi_{j_0\dots j_{n-1}})^{-1}\cdot\tau[x_0\dots x_{l-1}]\cdot
\tr (D\psi_{j_0\dots j_{n-1}})^{-1}  .  $$
This proves (2.8) when $x_0x_1\dots x_{n-1}=j_0j_1\dots j_{n-1}$ and we are done. 

\fin

We begin to note that (2.2) with $Q=Id$ easily implies that  
$$\kappa=(Id,\tau)_{HS}\le ||\tau||  .  \eqno (2.11)$$
Together with (2.6) this implies that $||\tau||$ and $\kappa$ are mutually absolutely continuous; since the polar decomposition of $\tau$ is 
$$\tau=T_x||\tau||$$
with $||T_x||_{HS}=1$ for $||\tau||$-a. e. $x\in S$, we get that  
$$\tau=\l_x T_x\kappa$$
where $\fun{\l}{S}{(0,+\infty)}$ is the Radon-Nikodym derivative 
$\frac{\dr||\tau||}{\dr\kappa}$; it is bounded away from zero and $+\infty$ by (2.6) and (2.11). 

By [13] we have that there is a Borel vector field $\fun{v}{S}{\R^2}$ of unitary vectors such that 
$$T_x=P_{v(x)}\txt{for $||\tau||$-a. e. $x\in S$}      $$
where $P_v$ is the projection we defined in the introduction. By the last two formulas, 
$$\tau=\l_x P_{v(x)}\kappa   .    $$
By Rokhlin's theorem (see [17], [18] or theorem 0.5.1 of [14]) the first equality below holds for $\kappa$-a. e. 
$x=(x_0x_1\dots)\in S$; the second equality follows from point 5) of proposition 2.1 and the definition of $\kappa$. 
$$\l_x P_{v(x)}=
\lim_{l\tends+\infty}
\frac{\tau[x_0\dots x_{l-1}]}{\kappa[x_0\dots x_{l-1}]}=$$
$$\lim_{l\tends+\infty}
\frac{
(D\psi_{x_0\dots x_{l-1}})\cdot\tau(S)\cdot\tr (D\psi_{x_0\dots x_{l-1}})
}{
tr[
(D\psi_{x_0\dots x_{l-1}})\cdot\tau(S)\cdot\tr (D\psi_{x_0\dots x_{l-1}})
]
}    .  \eqno (2.12)$$
As we recalled in proposition 2.1, $\tau(S)=Id$; thus, the last formula reduces to 
$$\l_x P_{v(x)}=\lim_{l\tends+\infty}
\frac{
(D\psi_{x_0\dots x_{l-1}})\cdot\tr (D\psi_{x_0\dots x_{l-1}})
}{
||
(D\psi_{x_0\dots x_{l-1}})
||^2_{HS}
}   .   \eqno (2.13)$$
In particular, the Hilbert-Schmidt norms on both sides coincide, implying that $\l_x\equiv 1$. 

We define the map 
$$\fun{H}{S\setminus \{ a,b,c \}}{\L(\R^d)}$$
$$H(x)= D\psi_i(x)\txt{if}x\in\psi_i(S)\setminus\{ a,b,c \}  .  $$
This function $H$ is well defined, since all the points which belong to more than one 
$\psi_i(S)$ are in $\{ a,b,c \}$; it is also defined $\kappa$-a. e. since we saw above that 
$\kappa(\{ a,b,c \})=0$. 

We set 
$$H^l(x)=H(x)\circ H(F(x))\circ\dots\circ H(F^{l-1}(x))       \eqno (2.14)$$
and note that $H^l$ is defined save on the countable set 
$$\bigcup_{l\ge 0}\bigcup_{x_0\dots x_{l-1}}\psi_{x_0\dots x_{l-1}}(\{ a,b,c \})    $$
i. e. save on the points where the coding $\Phi$ is not injective. 

If $x=\Phi(w_0w_1\dots)$, using (1.9) we see that, for $x$ in the full measure set where $H^l$ is defined,  
$$H^l(x)=D\psi_{w_0\dots w_{l-1}}  .  \eqno (2.15)$$

For the following lemma we recall that, if $v,w,a\in\R^d$ and $||v||=1$, then 
$$(v\otimes w)a=(a,v)\cdot w   .  \eqno (2.16)$$

\lem{2.3} Let the vector field $v(x)$ be as in (2.13). Then, for each $l\ge 1$ there is a vector field of unit vectors $v^l(x)$ such that, for $\kappa$-a. e. $x\in S$, 
$$||\frac{H^l(x)}{||H^l(x)||_{HS}}-v^l(x)\otimes v(x)||_{HS}\tends 0  .  \eqno (2.17)$$
Equivalently, if $x=\Phi(w_0w_1\dots)$, we get by (2.15) that 
$$||
\frac{D\psi_{w_0\dots w_l}}{||D\psi_{w_0\dots w_l}||_{HS}}-
v^l(x)\otimes v(x)
||_{HS}\tends 0  .  \eqno (2.18)$$
We are not asserting that the convergence is uniform. 

\proof Let $x=\Phi(w_0w_1\dots)$; the first equality below is the definition of $v^l$, the second one follows from (2.15). 
$$v^l(x)\colon=\frac{1}{||\tr H^l(x)||_{HS}}\tr H^l(x) v(x)=
\frac{
\tr(D\psi_{w_0\dots w_{l-1}})v(x)
}{
||\tr(D\psi_{w_0\dots w_{l-1}})||_{HS}
}  .  $$
By (2.13) and the fact that $\l_x\equiv 1$ we see that, for $\kappa$-a. e. $x\in G$, 
$$(v^l(x),v^l(x))\tends(P_{v(x)}v(x),v(x))=1  .  \eqno (2.19)$$
If we take $w\perp v(x)$ we obtain the equality below, while the limit follows from (2.13) and (2.15). 
$$\left(
\frac{\tr H^l(x)w}{||\tr H^l(x)||_{HS}},\frac{\tr H^l(x)w}{||\tr H^l(x)||_{HS}}
\right)  \tends
(P_{v(x)}w,w)=0  .  $$
In other words, for $l$ large $\frac{\tr H^l(x)}{||\tr H^l(x)||_{HS}}$ brings the unit vector 
$v(x)$ into the unit vector $v^l(x)$, and annihilates $v(x)^\perp$; this is tantamount to saying that 
$$||\frac{\tr H^l(x)}{||\tr H^l(x)||_{HS}}- v(x)\otimes v^l(x) ||_{HS}\tends 0 $$
for $\kappa$-a. e. $x\in S$. Transposing, we get (2.17). 

\fin

\vskip 2pc
\centerline{\bf \S 3}
\centerline{\bf The cells}
\vskip 1pc

Let $x=\Phi(w_0w_1\dots)\in S$ and let the cell $[w_0\dots w_l]$ be defined as in (1.11). Let $\bar A$, $\bar B$ and $\bar C$ be the three points of the "boundary" of $S$ as in section 1; we shall say that 
$$\bar A_{w,l}=\psi_{w,l}(\bar A),\qquad
\bar B_{w,l}=\psi_{w,l}(\bar B),\qquad
\bar B_{w,l}=\psi_{w,l}(\bar B)   \eqno (3.1)$$
are the three points of the boundary of $[w_0\dots w_l]$. 
Let us denote by $\tilde S$ the set of the points of $S$ where (2.13) holds; we saw before formula (2.12) that $\tilde S$ is a Borel set and $\kappa(\tilde S)=1$. Let $\hat S$ denote $\tilde S$ minus the countable set of points of (3.1); again $\hat S$ is Borel; since 
$\kappa$ is non-atomic, $\kappa(\hat S)=1$.

In section 1 we saw that $S$ is contained inside the closed set $M$ of (1.12); we define by $d(x,A)$ the Euclidean distance of $x\in\R^2$ from a set $A\subset\R^2$ and we set 
$$S_\theta=\{
x\in S\st d(x,M)\ge\theta
\}  .  \eqno (3.2)$$
Heuristically, $S_\theta$ is $S$ minus a thin strip around its three "sides", i . e. the curve $M$.

The definition of $S_\theta$ in (3.2) easily implies that 
$$\bigcup_{n\ge 1} S_\frac{1}{n}=S\setminus M  .  $$
By [8], $\kappa(M)=0$; together with the formula above this implies that there is 
$\theta_0>0$ such that, for $\theta\in(0,\theta_0]$, 
$S_\theta$ is not empty. 

Again since $\kappa(M)=0$, (2.6) implies that $||\tau||(M)=0$; as a consequence, 
$$\tau(S_\theta)\tends \tau(S)=Id
\txt{as}\theta\tends 0  .  \eqno (3.3)$$
As a last bit of notation, we shall denote by $P_v^\perp$ the projection on the orthogonal subspace to $v$: 
$$P_v^\perp w=w-P_v(w)  .  $$

\lem{3.1} Let $w=(w_0w_1\dots)\in\Sigma$ be such $x=\Phi(w)\in\hat S$; let $\fun{v}{\hat S}{\R^2}$ be the vector field of (2.13). Then, we can label the boundary points $A_{w,n}$, $B_{w,n}$ and $C_{w,n}$ of (3.1) in such a way that the two assertions below hold. 

\noindent 1) For $n\tends+\infty$ we have that
$$\frac{
\max[
||P_{v(x)}^\perp(A_{w,n}-B_{w,n})||, \quad
||P_{v(x)}^\perp(B_{w,n}-C_{w,n})||, \quad
||P_{v(x)}^\perp(C_{w,n}-A_{w,n})||
]
}{
||P_{v(x)}(A_{w,n}-B_{w,n})||
}  \tends 0   .  \eqno (3.4)$$
Heuristically, the cell $\psi_{w,n}(S)$ is "long" in the direction $v(x)$ and "short" in the orthogonal direction; the vertices are re-labeled in such a way that the side 
$A_{w,n}B_{w,n}$ has the longest projection on $v(x)$. 

\noindent 2) Let $\theta\in(0,\theta_0]$, let $S_\theta$ be as in (3.2) and let us suppose that $x\in \psi_{w,n_k}(S_\theta)$ for a sequence $n_k\nearrow+\infty$; then, 
$$\max\left(
\frac{
||P_{v(x)}^\perp(x-A_{w,n_k})||
}{
||P_{v(x)}(x-A_{w,n_k})||
},\quad
\frac{
||P_{v(x)}^\perp(x-B_{w,n_k})||
}{
||P_{v(x)}(x-B_{w,n_k})||}
\right)
\tends 0  .   \eqno (3.5)  $$
Again heuristically, $x$ sees both boundary points on the "long" side of the cell much farther away in the $v(x)$ direction than in the orthogonal one. Note that neither (3.4) nor (3.5) are uniform in $x\in\hat S$. 

\proof We begin with (3.4). Let $\bar A\bar B\bar C$ be the same equilateral triangle of section 1; its sides are $\frac{2}{\sqrt 3}$ long, its height $1$. Elementary geometry tells us that, given any 
$v\in\R^2\setminus\{ 0 \}$, 
$$\max(
||P_v(\bar B-\bar A)||, \quad ||P_v(\bar B-\bar C)||, \quad ||P_v(\bar C-\bar A)||
)=
\diam(P_v(\bar A\bar B\bar C))\ge 1  .  $$
Clearly, the diameter of the projection is minimal when $v$ is parallel to one of the heights of $\bar A\bar B\bar C$. 

Let $x=\Phi(w)\in\hat S$ and let  $v^n(x)$ and $v(x)$ be as in lemma 2.3; we can relabel 
$\bar A$, $\bar B$, $\bar C$ to $A$, $B$, $C$ in such a way that $||P_{v^n(x)}(B-A)||$ is maximal; using (2.16), and the fact that $||v(x)||=1$ we get the equality below; the first inequality comes from the formula above; the second on comes, for some 
$\d_n\searrow 0$, from (2.19).

$$||v^n(x)\otimes v(x)(B-A)||= (v^n(x),B-A)\ge$$
$$||v^n(x)|| \ge 1-\d_n  .  \eqno (3.6)  $$
As in (3.1), we define 
$$A_{w,n}=\psi_{w,n}(A),\qquad
B_{w,n}=\psi_{w,n}(B),\qquad
C_{w,n}=\psi_{w,n}(C)  .  $$
Together with the fact that $\psi_{w,n}$ is affine, this definition implies the equality below.  Since $x\in\hat S$, (2.18) holds; this implies that, for $n\ge n_0(x)$, the first inequality follows for a matrix $F_n$ with $||F_n||_{HS}\le\e$. This bound on the norm implies the second inequality, while the last one comes from (3.6). 
$$\frac{
||P_{v(x)}(A_{w,n}-B_{w,n})||
}{
||D\psi_{w,n}||_{HS}
}  =
\frac{
||P_{v(x)}\circ D\psi_{w,n}(A-B)||
}{
||D\psi_{w,n}||_{HS}
} \ge  $$
$$||v^n(x)\otimes v(x)(A-B)||-||P_{v(x)}\circ F_n(A-B)||\ge$$
$$||v^n(x)\otimes v(x)(A-B)||-\e||A-B||\ge 1-\d_n-\frac{2}{\sqrt 3}\e  .  \eqno (3.7)$$
The same argument shows that, for $n$ large, the formula below holds. 
$$\frac{
||P_{v(x)}^\perp(A_{w,n}-B_{w,n})||
}{
||D\psi_{w,n}||_{HS}
}   =  
\frac{||P_{v(x)}^\perp\circ D\psi_{w,n}(A-B)||
}{
||D\psi_{w,n}||_{HS}
}  = $$
$$||P_{v(x)}^\perp F_n(A-B)||\le
\e  ||A-B||=\e  \frac{2}{\sqrt 3}  .  \eqno (3.8)$$
Clearly, the same inequality holds for the other two couples of vertices; now (3.4) follows from (3.7), (3.8) and the fact that $\d_n\tends 0$. 

We prove (3.5). For starters, we assert that there is $a(\theta)>0$ such that, if 
$\tilde x\in S_\theta$, 
$$||P_{v^n(x)}(\tilde x-A)||  \ge a(\theta)\cdot||\tilde x-A||  \txt{and}
||P_{v^n(x)}(\tilde x-B)||  \ge a(\theta)\cdot||\tilde x-B||  .    \eqno (3.9)$$
This follows by our choice of $A$ and $B$ in (3.6). Indeed, since $S_\theta$ is a compact set contained in the open, equilateral triangle $ABC$, we see that there is $b(\theta)>0$ such that the angle between $AB$ and $A\tilde x$ belongs to 
$[b(\theta),\frac{\pi}{3}-b(\theta)]$. Since $P_{v^n(x)}=P_{-v^n(x)}$, we can as well suppose that $(v^n(x),B-A)>0$; since $AB$ is the side with the longest projection on 
$v^n(x)$, the angle between $AB$ and $v^n(x)$ lies in $[-\frac{\pi}{6},\frac{\pi}{6}]$. Summing up, we get that that the angle between $v^n(x)$ and $A\tilde x$ lies in 
$[-\frac{\pi}{6}+b(\theta),\frac{\pi}{2}-b(\theta)]$. The first formula of (3.9) now follows with $a(\theta)=\sin(b(\theta))$; the second one follows analogously for the same $a(\theta)$. 

Since $x\in\psi_{w,n}(S_\theta)$, there is $\tilde x\in S_\theta$ such that 
$x=\psi_{w,n}(\tilde x)$. The first inequality below follows as in (3.7) for a matrix $F_n$ with $||F_n||_{HS}\le\e$, this latter fact implying the second inequality; the last one follows, for some $\d_n\searrow 0$, from (3.9) and (2.19). 
$$\frac{
||P_{v(x)}(x-A_{w,n})||
}{
||D\psi_{w,n}||_{HS}
}  \ge 
||{v^n(x)}\otimes v(x) (\tilde x-A)||-||P_{v(x)}\circ F_n(\tilde x-A)||\ge$$
$$(v^n(x),\tilde x-A)-\e||\tilde x-A||\ge
[a(\theta)(1-\d_n)-\e]\cdot ||\tilde x-A||  .  $$
Analogously, 
$$\frac{
||P_{v(x)}(x-B_{w,n})||
}{
||D\psi_{w,n}||_{HS}
}  \ge
[a(\theta)(1-\d_n)-\e] \cdot ||\tilde x-B||  .  $$
Given $\e>0$, we see as in (3.8) that, for $n$ large enough (depending on $w$ and thus on $x$) we have 
$$\frac{
||P_{v(x)}^\perp(x-A_{w,n})||
}{
||D\psi_{w,n}||_{HS}
}  \le\e
||\tilde x-A||   ,   \qquad
\frac{
||P_{v(x)}^\perp(x-B_{w,n})||
}{
||D\psi_{w,n}||_{HS}
}  \le\e
||\tilde x-B||
  .  $$
Now (3.5) follows from the last three formulas. 

\fin

\lem{3.2} There is $\d(\theta)\tends 0$ as $\theta\tends 0$ such that, for all 
$x=\Phi(w)\in\hat S$, 
$$\lim_{n\tends+\infty}
\frac{
\kappa(\psi_{w,n}(S_\theta))
}{
\kappa(\psi_{w,n}(S))
}  \ge
1-\d(\theta)  .  \eqno (3.10)$$
As in lemma 3.1, we are {\rm not} asserting that the limit is uniform in $w\in\Sigma$. 

\proof The first equality below comes from lemma 2.2. The second equality comes from the definition of the push-forward.  
$${\b^n}\cdot
\frac{1}{
|| D\psi_{w,n} ||_{HS}^2
}  \cdot \tau(\psi_{w,n}(S_\theta))  =  $$
$$\frac{
D\psi_{w,n}
}{
|| D\psi_{w,n} ||_{HS}
}  \cdot
[(\psi_{w,n})_\sharp\tau](\psi_{w,n}(S_\theta))\cdot
\frac{
\tr D\psi_{w,n}
}{
|| D\psi_{w,n} ||_{HS}
}  =$$
$$\frac{
D\psi_{w,n}
}{
|| D\psi_{w,n} ||_{HS}
}  \cdot
\tau(S_\theta)\cdot
\frac{
\tr D\psi_{w,n}
}{
|| D\psi_{w,n} ||_{HS}
}  .  $$
If $A$ lies in a bounded set of $M^d$, the family of maps $\fun{\a_A}{B}{BA\tr B}$ is equicontinuous on the set $||B||_{HS}\le 1$. As a consequence, the maps $\fun{\a_\theta}{B}{B\tau(S_\theta)\tr B}$ are equicontinuous as $\theta$ varies in $[0,1]$. Together with (2.18) this implies the limit below, while the equality comes from the last formula.  
$$\sup_{\theta\in[0,1]}
||{\b^n}\cdot
\frac{1}{
|| D\psi_{w,n} ||_{HS}^2
}  \cdot \tau(\psi_{w,n}(S_\theta)) -
v^n(x)\otimes v(x)\cdot\tau(S_\theta)\cdot\tr (v^n(x)\otimes v(x))||_{HS}=$$
$$\sup_{\theta\in[0,1]}
||\frac{
D\psi_{w,n}
}{
|| D\psi_{w,n} ||_{HS}
}  \cdot
\tau(S_\theta)\cdot
\frac{
\tr D\psi_{w,n}
}{
|| D\psi_{w,n} ||_{HS}
} -
v^n(x)\otimes v(x)\cdot\tau(S_\theta)\cdot\tr (v^n(x)\otimes v(x))||_{HS}
\tends 0    
\eqno (3.11)$$
as $n\tends+\infty$. 

We recall that the map $\fun{\r_B}{A}{\tr BAB}$ is Lipschitz with a Lipschitz constant which only depends on $||B||_{HS}$; next we note that $||v^n(x)\otimes v(x)||$ is bounded, since $v^n(x)$ is bounded by definition; by this and (3.3) there is 
$\d^\pprime(\theta)\tends 0$ as $\theta\tends 0$, independent of $n$, such that 
$$||
v^n(x)\otimes v(x)\cdot\tau(S_\theta)\cdot v(x)\otimes v^n(x)-
v^n(x)\otimes v(x)\cdot\tau(S)\cdot v(x)\otimes v^n(x)||_{HS}\le \d^\pprime(\theta)  .  $$
Since $\tau(S)=Id$ by proposition 2.1, this is tantamount to 
$$||
v^n(x)\otimes v(x)\cdot\tau(S_\theta)\cdot v(x)\otimes v^n(x)-
P_{v^n(x)}||_{HS}\le \d^\pprime(\theta)  .  $$
By the last formula and (3.11) we get that, for $n$ large enough (depending on 
$x\in\hat S$),
$$||
\b^n\cdot\frac{1}{||D\psi_{w,n}||^2_{HS}}\cdot\tau(\psi_{w,n}(S_\theta))-
P_{v^n(x)}
||_{HS}  \le 2\d\pprime(\theta)  .  $$
Taking the Hilbert-Schmidt product with the identity we get that  
$$\left\vert
\left(
{\b^n}\cdot
\frac{1}{
|| D\psi_{w,n} ||_{HS}^2
}  \cdot \tau(\psi_{w,n}(S_\theta)) , Id
\right)_{HS}  -1
\right\vert  <2\d^\pprime(\theta)  .  $$
Analogously, 
$$\left\vert
\left(
{\b^n}\cdot
\frac{1}{
|| D\psi_{w,n} ||_{HS}^2
}  \cdot \tau(\psi_{w,n}(S)),Id
\right)_{HS}  -1
\right\vert\le 2\d^\pprime(\theta)  .  $$
The first equality below comes from the definition of Kusuoka's measure $\kappa$; the second one comes since we are integrating the constant matrix field $Id$; the inequality comes from the last two formulas; it holds for $n$ large enough, depending on $x\in\hat S$. 
$$\left\vert
\frac{\kappa(\psi_{w,n}(S_\theta))}{\kappa(\psi_{w,n}(S))}-1
\right\vert=
\left\vert
\frac{
\int_{\psi_{w,n}(S_\theta)}(Id,\dr\tau(x))_{HS}
}{
\int_{\psi_{w,n}(S)}(Id,\dr\tau(x))_{HS}
}  -1
\right\vert
=     $$
$$\left\vert
\frac{
(Id,\tau(\psi_{w,n}(S_\theta)))_{HS}
}{
(Id,\tau(\psi_{w,n}(S)))_{HS}
}-1
\right\vert   \le 4\d^\pprime(\theta)  .  $$
If we set $\d(\theta)\colon=4\d^\pprime(\theta)$,  (3.10) follows from the last formula.

\fin

\vskip 1pc

\noindent{\bf Definition.} Let $\hat S$ be defined as at the beginning of this section and let $\d(\theta)>0$ be as in lemma 3.2. We define $F_\theta$ as the set of all the 
$x=\Phi(w)\in\hat S$ such that for infinitely many $n$'s we have 
$$x\in\psi_{w,n}(S_\theta)
\txt{and}
\frac{
\kappa(\psi_{w,n}(S_\theta))
}{
\kappa(\psi_{w,n}(S))
}  \ge 1-2\d(\theta)   .  \eqno (3.12)$$
Recalling the definition of $S_\theta$ in (3.2), the first formula above says that $x$ is "well in the interior" of the cell $\psi_{w,n}(S)$. The second inequality holds for all 
$x=\Phi(w)\in\hat S$ and $n\ge n_0(w)$ by lemma 3.2.

\lem{3.3} Let $\theta_0$ be defined as after (3.2) and let the function $\d(\theta)$ be as in lemma 3.2; then,  
$$\kappa(F_\theta)\ge 1-2\d(\theta)
\txt{for all $\theta\in(0,\theta_0]$.}     \eqno (3.13)$$

\proof Let us consider the set of the $x=\Phi(w)\in\hat S$ which satisfy (3.12) for some fixed
$n\in\N$; clearly, it is  
$$F_{n,\theta}\colon=
\bigcup_{w_0\dots w_{n}}\psi_{w_0\dots w_{n}}(S_\theta)  \cap\hat S  $$
where the union is over all words $w_0\dots w_{n}$ for which the second formula of (3.12) holds. Note that the union is disjoint: indeed, the two cells $\psi_{w_0\dots w_n}(S)$ and 
$\psi_{w_0^\prime\dots w_n^\prime}(S)$ either coincide (and then 
${w_0\dots w_n}={w_0^\prime\dots w_n^\prime}$) or are disjoint or intersect at the boundary points; but $\psi_{w_0^\prime\dots w_n^\prime}(S_\theta)$ does not contain the boundary points.

Recall that, by definition, 
$$F_\theta=\limsup_{n\tends+\infty} F_{n,\theta}  .  $$
Thus, it suffices to show that  
$$\kappa(F_{n,\theta})\ge 1-2\d(\theta)  $$
from a certain $n$ onwards; equivalently, we are going to show that 
$$\kappa(F_{n,\theta}^c)\le 2\d(\theta)  .  \eqno (3.14)$$

We saw above that the cells $\psi_{w_0\dots w_n}(S)\cap\hat S$  are disjoint. In other words, if $x\in\hat S$, $x$ lies in only one $\psi_{w_0^\prime\dots w_n^\prime}(S)$; thus, 
$x\not\in F_{n,\theta}$ only in two cases:

\noindent 1) $x\in\psi_{w_0\dots w_n}(S)$ with the word $(w_0\dots w_n)$ which does not satisfy the second formula of (3.12), or

\noindent 2) $(w_0,\dots,w_n)$ satisfies the second formula of (3.12), but 
$x\in\psi_{w_0\dots w_n}(S\setminus S_\theta)$. 

By (3.10) and Rokhlin's theorem ([17], [18]) we see easily that the measure of the union of the cells 1) tends to zero as $n\tends+\infty$. Thus, it suffices to prove that the measure of the union of 
$\psi_{w_0\dots w_n}(S\setminus S_\theta)$ for the words satisfying 2) is smaller than 
$2\d(\theta)$.  Since these cells are disjoint (recall that the boundary points where they intersect are not in $\hat S$), we get the equality below, where the sum is over the words 
$w_0\dots w_{n}$ which satisfy the second formula of (3.12); the first inequality follows by (3.12), the second one from the fact that the union is disjoint and $\kappa$ is probability. 
$$\kappa\left(
\bigcup\psi_{w_0\dots w_{n}}(S\setminus S_\theta)
\right)   =
\sum\kappa(\psi_{w_0\dots w_{n}}(S\setminus S_\theta))  \le$$
$$2\d(\theta)\sum\kappa(\psi_{w_0\dots w_{n}}(S))\le
2\d(\theta)  .  $$
This is (3.14), which we have seen to imply the thesis.

\fin

\vskip 1pc

\centerline{\bf \S 4}

\centerline{\bf Cheeger's energy}

\vskip 1pc

We begin the proof of theorem 1 showing that, if $f\in C^1(\R^2,\R)$, then 
$Ch(f)\le 2\ec(f)$. The core of the proof is in lemma 4.3 below, which says the following: if 
$\pi$ is a test plan with bounded deformation then, for $\pi$-a. e. curve $\g$ and $\L^1$-a. e. $t\in[0,1]$, $\g$ is differentiable at $t$ and $\dot\g_t$ is parallel to the unit  vector $v(\g_t)$ of (2.13). Thus, if $\phi\in C^1(\R^2,\R)$, we get that $|(\nabla\phi(x),v(x))|$ is a weak gradient of $\phi$; by the definition of Cheeger's derivative, this implies that 
$|D\phi|_w(x)\le |(\nabla \phi(x),v(x))|$.

We set 
$$\hat F=\bigcup_{n\ge 1} F_\frac{1}{n}  \eqno (4.1)$$
where $F_\theta$ is defined as in formula (3.12). By lemma 3.3 we have that 
$\kappa(\hat F)=1$; since $F_\frac{1}{n}\subset\hat S$ by definition, we have that 
$\hat F\subset\hat S$; recall that $\hat S$ has been defined at the beginning of section 3.

\lem{4.1} Let the set $\hat F$ be defined as in (4.1); let $x=\Phi(w)\in\hat F$ and let 
$\g\in AC^2([0,1],S)$, $t\in(0,1)$ be such that $\g_t=x$ and $\g$ is differentiable at $t$. Let the vector field $v(x)$ be as in (2.13). Then, 
$$\dot\g_t=\pm ||\dot\g_t||\cdot v(\g(t))  .  \eqno (4.2)$$

\proof By the definition of $\hat F$ in (4.1) we have that $x\in F_\frac{1}{l}$ for some 
$l\in\N$. By the definition of $F_\frac{1}{l}$ in (3.12) there is a sequence 
$n_k\nearrow+\infty$ such that 
$$x\in\psi_{w,n_k}(S_\frac{1}{l})   \qquad\forall k\ge 1  .  \eqno (4.3)$$

Since $x\in\hat S$, we can label the boundary of the cell $\psi_{w,n_k}(S)$ as in lemma 3.1; we get that (3.5) holds for the point  $x=\Phi(w)$. Since $\g_t=x$ and $\g$ can escape out of a cell only hitting its boundary, three cases are possible.

\noindent 1) There is $h_k\tends 0$ such that 
$$\g(t+h_k)\in\{ A_{w,n_k},B_{w,n_k} \}  .  $$
By (4.3), (3.5) holds and implies that 
$$\frac{\g(t+h_k)-\g(t)}{h_k}\tends\l v(x)  $$
for some $\l\in\R$; since we are supposing that $\g$ is differentiable at $t$, (4.2) follows. 

\noindent 2) There are $h_k^+>0$, $h_k^-<0$ such that $h_k^+,h_k^-\tends 0$ and 
$$\g(t+h_k^+)=\g(t+h_k^-)=C_{w,k}  .  $$
Since $\g$ is differentiable at $t$, we get the first equality below, while the second one follows by the formula above.
$$\dot\g_t=
\lim_{k\tends+\infty}
\frac{\g(t+h_k^+)-\g(t+h_k^-)}{h_k^+-h_k^-}=0  .  $$
This implies (4.2) for the speed $||\dot\g_t||=0$. 

\noindent 3) The only remaining alternative is that, for some $h_0>0$, either 
$$\g([t,t+h_0])\subset\psi_{w,n_k}(S)
\txt{or}
\g([t-h_0,t])\subset\psi_{w,n_k}(S)  $$
for infinitely many $k$'s. 

Let us suppose that the first case happens and let us fix $h\in(0,h_0]$. Since the diameter of $\psi_{w,n_k}(S)$ tends to zero by (1.3), the limit below holds, while the inequality holds for $h\in(0,h_0]$.
$$\frac{
||\g(t+h)-\g(t)||
}{
h
}\le
\frac{
\diam(\psi_{w,n_k}(S))
}{
h
}  \tends 0  .  $$
Thus, 
$$\frac{||\g(t+h)-\g(t)||}{h}=0$$
for all $h\in(0,h_0]$. Letting $h\tends 0$ and recalling that $\g$ is differentiable at $t$, we get that $\dot\g_t=0$, proving (4.2) also in this case. 

\fin

\lem{4.2} Let $\pi$ be a test plan with bounded deformation and let $\hat F\subset S$ be defined as (4.1). Then, denoting by $\L^1$ the Lebesgue measure on $\R$, we have that 
$$\L^1(\g^{-1}(\hat F^c))=0
\txt{for $\pi$-a. e. $\g\in C([0,1],S)$.}  \eqno (4.4)$$

\proof We consider $[0,1]\times C([0,1], S)$ with the measure $\L^1\otimes\pi$; we define the map 
$$\fun{g}{[0,1]\times C([0,1], S)}{[0,1]\times S}  $$
$$\fun{g}{(t,\g)}{(t,\g_t)}  .  $$
Since $\pi$ has bounded deformation, it is immediate that there is $C>0$ such that 
$$g_\sharp(\L^1\otimes\pi)\le C(\L^1\otimes \kappa)  .  $$
We saw after (4.1) that $\kappa(\hat F^c)=0$, which implies by Fubini that  $(\L^1\otimes \kappa)([0,1]\times\hat F^c)=0$. Together with the last formula this implies that $g_\sharp(\L^1\otimes \pi)([0,1]\times\hat F^c)=0$, i. e., by the definition of push-forward, 
$$(\L^1\otimes\pi)(\{ (t,\g)\st\g_t\not\in\hat F \})  =0  .  $$
By Fubini, this implies the thesis.

\fin

\lem{4.3} Let $\phi\in C^1(\R^2,\R)$. Then, 
$$|D\phi|_w(x)\le |(P_{v(x)},\nabla \phi(x))|
\txt{for $\kappa$-a. e. $x\in S$.}  \eqno (4.5)$$

\proof
By the definition of $|D\phi|_w(x)$ it suffices to show that $|(P_{v(x)},\nabla \phi(x))|$ is a weak gradient of $\phi$; by the definition of weak gradient in section 1, this follows if we show that, for all test plans $\pi$ with bounded deformation, for $\pi$-a.e. $\g\in C([0,1],S)$ and for $\L^1$-a. e. $t\in[0,1]$ which is a point of differentiability of $\g$ we have 
$$\dot\g_t=\pm ||\dot\g_t||\cdot v(\g(t))  .  \eqno (4.6)$$

Let $\pi$ be a test plan and let us call $B\subset C([-1,1],S)$ the full measure set of curves for which (4.4) holds; since $\pi$ is a test plan, excluding a zero-measure set we can suppose that all curves in $B$ are absolutely continuous. If we set $A_\g=\g^{-1}(\hat F)$, we see by (4.4) that 
$\L^1(A_\g)=1$ for all $\g\in B$. Since $\g$ is absolutely continuous, if we set 
$$\hat A_\g=\{  \quad
t\in A_\g\st\g\txt{is differentiable at $t$}
\}  $$
we see that $\L^1(\hat A_\g)=1$ too. The upshot of all this is that it suffices to show (4.6) for all $\g\in B$ and all $t\in\hat A_\g$. To say that $t\in\hat A_\g$ is equivalent to say that 
$\g_t\in\hat F$ and $\g$ is differentiable at $t$; now (4.6) follows from lemma 4.1. 

\fin

\vskip 1pc

We need to prove the opposite inequality to (4.5); for this we shall use lemma 3.1 and Cheeger's characterisation of the energy (1.13). We begin with two lemmas; before stating them we recall from the end of section 1 that 
$Lip_a(f,x)$ is the local Lipschitz constant {\sl calculated on $S$}; the distance on $S$ is that induced by the immersion in $\R^2$. 

\lem{4.4} Let $f\in Lip(S)$. Then, 
$$Ch(f)=\inf\liminf_{n\tends+\infty}\2\int_SLip_a^2(f_n)\dr\kappa  $$
where the $\inf$ is over all sequences $\{ f_n \}_{n\ge 1}$ such that 

\noindent 1) $f_n\in C^1(\R^2,\R)$ and 

\noindent 2) $f_n\tends f$ in $L^2(S,\kappa)$. 

\proof Since $C^1$ functions are Lipschitz, the inequality 
$$Ch(f)\le\inf\liminf_{n\tends+\infty}\2\int_SLip_a^2(f_n)\dr\kappa  $$
follows immediately from (1.13). In the next steps we show the opposite one. 

\noindent {\bf Step 1.} We begin to assert that in (1.13) it suffices to consider sequences 
$\{ g_n \}_{n\ge 1}$ of Lipschitz functions which are $C^1$ on each cell 
$\psi_{w_1\dots w_n}(\bar A\bar B\bar C)$. 

We prove the assertion. Let $\fun{h}{S}{\R}$ be Lipschitz and let us consider the cell 
$\psi_{w_0\dots w_n}(S)$; its three boundary points $A_{w_0\dots w_n}$, 
$B_{w_0\dots w_n}$ and $C_{w_0\dots w_n}$ define a triangle. It is standard (Kirszbraun's theorem, see for instance [6]) that we can extend $h$ restricted to the three points 
$A_{w_0\dots w_n}, B_{w_0\dots w_n}, C_{w_0\dots w_n}$  to a function $\tilde h$ defined on the whole triangle $A_{w_0\dots w_n}, B_{w_0\dots w_n}, C_{w_0\dots w_n}$ which has the same Lipschitz constant $h$ has on the three points. By definition, $h$ and 
$\tilde h$ coincide on the three boundary points. 

We use another standard approximation: there exists $\e_n\tends 0$ and a function $g_n$ such that $g_n$ is $C^1$ in the interior of each triangle $\psi_{w_0\dots w_n}(S)$ and satisfies 

\noindent $i$) $g_n(A_{w_0\dots w_n})= h(A_{w_0\dots w_n})$, 
$g_n(B_{w_0\dots w_n})= h(B_{w_0\dots w_n})$ and 
$g_n(C_{w_0\dots w_n})= h(C_{w_0\dots w_n})$. 

\noindent $ii$) $|g_n(x)-\tilde h(x)|\le\e_n$ for all $x$ in the triangle.

\noindent $iii$) $|\nabla g_n(x)|\le L_{w_0\dots w_n}+\e_n$ for all $x$ in the triangle, where $L_{w_0\dots w_n}$ is the Lipschitz constant of $\tilde h$ in the triangle (or of $h$ on the three points $A_{w_0\dots w_n}, B_{w_0\dots w_n}, C_{w_0\dots w_n}$, which is the same.) Here and in the following we denote by $\nabla f$ the standard gradient of $f$ in 
$\R^2$. 

For a fixed $n\in\N$, the triangles $\psi_{w_0\dots w_n}(ABC)$ intersect only at the vertices; this implies by $i$) that, if we define $\tilde g_n$ to be equal to $g_n$ on each triangle, the definition is well-posed and the function $\tilde g_n$ is continuous. Moreover, it is $C^1$ in the interior of each triangle. 

Since $h$ is Lipschitz, $ii$) and the definition of $\tilde h$ imply that $g_n\tends h$ uniformly on $S$, and thus in $L^2(S,\kappa)$. On the other side, by $iii$) and the definition of the $Lip_a$ in section 1 we have that 
$$\limsup_{n\tends+\infty}|\nabla g_n|(x)\le Lip_a(h,x)  $$
for all $x\in S$, save for the points on the boundary of a cell, where the gradient is not defined; but these are a countable set and have zero measure. 

Since $h$ is Lipschitz, we get that $L_{w_0\dots w_n}$ is bounded; now $iii$), the formula above and dominated convergence imply that 
$$\limsup_{n\tends+\infty}\int_S|\nabla g_n|^2(x)\dr\kappa(x)\le
\int_S Lip_a^2(h,x)\dr\kappa(x)  .  \eqno (4.7)$$
Now, let $\{ h_n \}_{n\ge 1}$ be a sequence of Lipschitz functions such that 
$$\2\int_S Lip_a^2(h_n,x)\dr\kappa(x)\tends Ch(f)      $$
and 
$$||h_n-f||_{L^2(S,\kappa)}\tends 0  .  $$
By the aforesaid we can approximate each $h_n$ with a function $g_n$ which is 

\noindent $\bullet$) continuous on $S$,

\noindent $\bullet$) $C^1$ on the $k$-th generation cells for some large $k=k(n)$ 

\noindent $\bullet$) and such that, by (4.7), 
$$\2\int_S|\nabla g_n|^2(x)\dr\kappa(x)\le
\2\int_S Lip_a^2(h_n,x)\dr\kappa(x)+\frac{1}{n}    $$
and 
$$||h_n-g_n||_{L^2(S,\kappa)}\tends 0  .  $$

From the last four formulas we get that $g_n\tends f$ in $L^2(S,\kappa)$ and that
$$\2\int_S|\nabla g_n|^2(x)\dr\kappa(x)\tends Ch(f)   .   $$
Clearly, this implies step 1. 

\noindent{\bf Step 2.} Next, we assert that in (1.13) we can take $f_n\in C^1(\R^2,\R)$. 

In order to show this, let us fix $n$ and let us gather all the boundary points of the $n$-th generation cells $[w_0\dots w_n]$ in a finite set $\{ P_j \}_{j=1}^{l_n}$. For $r>0$ we set 
$$O(r)=S\cap
\bigcup_{j=1}^{l_n} B(P_j,r_n)   .  $$
Since $\kappa$ is non-atomic, given $\e_n>0$ we can find $r_n>0$ such that 
$$\kappa(O(r_n))<\e_n  .  \eqno (4.8)$$
Let $g_n$ be the function of the previous step, which is $C^1$ on all the $n$-th generation triangles; it is easy to see that, whatever is $\e_n>0$ (or $r_n>0$, which depends on 
$\e_n$), we can find $f_n\in C^1(\R^2,\R)$ which coincides with $g_n$ on 
$S\setminus O(r_n)$ and such that, for all $x\in S$, $|\nabla f_n|(x)\le 2 Lip(g_n)$. Choosing $\e_n$ small enough in (4.8), this implies first that $f_n\tends f$ in $L^2(S,\kappa)$, second that 
$$\int_S||\nabla f_n|^2-|\nabla g_n|^2|\dr\kappa\tends 0$$
as $n\tends+\infty$. As at the end of step 1 this implies that 
$$\2\int_S|\nabla f_n|^2(x)\dr\kappa(x)\tends Ch(u)$$
which in turn implies (1.13) for $\{ f_n \}_{n\ge 1}$. 

\fin

\lem{4.5} Let $f\in C^1(\R^2,\R)$ and let the local Lipschitz constant $Lip_a(f,x)$ be defined as above; let $v(x)$ be as in (2.13). We assert that 
$$Lip_a(f,x)\ge |P_{v(x)}\nabla f(x)|\qquad\forall x\in\hat S  .  \eqno (4.9)$$

\proof Let $\psi_{w,n}(S)$  be the $n$-th generation cell which contains $x$; we relabel its boundary $A_{w,n}$, $B_{w,n}$ and $C_{w,n}$ in such a way that (3.4) holds. Since 
$A_{w,n}$ and $B_{w,n}$ belong to $S$ and converge to $x$ we have that 
$$Lip_a(f,x)\ge\lim_{n\tends+\infty} 
\frac{|f(B_{w,n})-f(A_{w,n})|}{||B_{w,n}-A_{w,n}||}  .  $$
Since $f\in C^1(\R^2,\R)$, (3.4) implies (4.9). 

\fin

\noindent{\bf End of the proof of theorem 1.} {\bf Step 1.} We assert that it suffices to show that $\ec(f)=2 Ch(f)$ when $f\in C^1(\R^2)$. 

We prove the assertion using theorem 3.8 of [8] (see section 4 of [11] for the original treatment and the proof). This theorem says that $\dc(\ec)$, the domain of $\ec$, is a Hilbert space for the inner product 
$$(u,v)_{\dc(\ec)}\colon=(u,v)_{L^2(S,\kappa)}+\ec(u,v)  $$
and that $C^1(\R^2)$ is dense in $\dc(\ec)$. 

Since $\ec$ and $2 Ch$ coincide on $C^1(\R^2)$, we get that $\dc(\ec)$ with the inner product $(\cdot,\cdot)_{\dc(\ec)}$ is also the completion of $C^1(\R^2)$ for the norm 
$$||u||^2_{\dc(Ch)}\colon=||u||^2_{L^2(S,\kappa)}+2Ch(u)  .  $$
This implies that $\dc(Ch)\supset\dc(\ec)$; the opposite inclusion follows by (1.13) if we show that $\dc(\ec)$ contains the Lipschitz functions. 

This follows from two facts: first, on $\dc(\ec)$, $\ec$ is defined by relaxation
$$\ec(u)=\inf\liminf_{n\tends+\infty}\ec(u_n)  .   \eqno (4.10)$$
The $\inf$ in the formula above is over all sequences 
$\{ u_n \}_{n\ge 1}\subset C^1(\R^2,\R)$ such that $u_n\tends u$ in $L^2(S,\kappa)$. 

The second fact is that, if $u\in Lip(S)$, then $u$ can be approximated, in the uniform topology, by a sequence $\{ u_n \}_{n\ge 1}\subset C^1(\R^2)$ such that 
$||\nabla u_n||_{\sup}$ is bounded. This implies the second inequality below; the first one comes from the fact that $\ec$ is lower semicontinuous; the equality comes from the definition of $\ec$ in proposition 2.1. 
$$\ec(u)\le\liminf_{n\tends+\infty}\ec(u_n)=
\liminf_{n\tends+\infty}\int_S(\nabla u_n(x),P_{v(x)}\nabla u_n(x))\dr\kappa (x)
<+\infty  .  $$

\noindent {\bf Step 2.} We show that $\ec(f)=2 Ch(f)$ when $f\in C^1(\R^2,\R)$. 

First of all, lemma 4.3 implies that 
$$2Ch(f)\le\ec(f)  .  \eqno (4.11)$$

We prove the opposite inequality. For $f\in C^1(\R^2,\R)$ , let us define its pre-Cheeger energy $pCh(f)$ as 
$$pCh(f)=\2\int_S Lip_a(f,x)^2\dr\kappa(x)  .  $$
We know by lemma 4.5 that 
$$pCh(f)\ge \2\ec(f)  .   $$
By lemma 4.4 we also have that, if $f\in C^1(\R^2,\R)$, 
$$Ch(f)=\inf\liminf_{n\tends+\infty}pCh(f_n)  $$
where the $\inf$ is over all sequences $\{ f_n \}_{n\ge 1}$ of $C^1$ functions such that 
$f_n\tends f$ in $L^2(S,\kappa)$. The last two formulas imply the first inequality below, while the second one comes from the fact that $\ec$ is lower semicontinuous. 
$$Ch(f)\ge
\inf\liminf_{n\tends+\infty}\2\ec(f_n)\ge
\2\ec(f)  .  $$

\fin

\vskip 2pc
\centerline{\bf References}





\noindent [1] L. Ambrosio, N. Gigli, G. Savar\'e, Gradient Flows, Birkh\"auser, Basel, 2005.

\noindent [2] L. Ambrosio, N. Gigli, G. Savar\'e, Heat flow and calculus on metric measure spaces with Ricci curvature bounded below - the compact case. Analysis and numerics of Partial Differential Equations, 63-115, Springer, Milano, 2013.  

\noindent [3] U. Bessi, Another point of view on Kusuoka's measure, preprint. 

\noindent [4] J. Cheeger, Differentiability of Lipschitz functions on metric measure spaces, GAFA, Geom. Funct. Anal. {\bf 9}, 428-517, 1999.

\noindent [5] S. Di Marino, Recent advances in BV and Sobolev spaces in metric measure spaces, Ph. D. thesis, 2014. 

\noindent [6] H. Federer, Geometric measure theory, Berlin, 1996. 



\noindent [7] A. Johansson, A. \"Oberg, M. Pollicott, Ergodic theory of Kusuoka's measures, J. Fractal Geom., {\bf 4}, 185-214, 2017.  

\noindent [8] N. Kajino, Analysis and geometry of the measurable Riemannian structure on the Sierpinski gasket, Contemporary Math. {\bf 600}, Amer. Math. Soc., Providence, RI, 2013.  

\noindent [9] J. Kigami, Analysis on fractals, Cambridge tracts in Math., {\bf 143}, Cambridge Univ. Press, Cambridge, 2001. 

\noindent [10] J. Kigami, Measurable Riemannian geometry on the Sierpinski gasket: the Kusuoka measure and the Gaussian heat kernel estimates, Math. Ann., {\bf 340}, 781-804, 2008. 

\noindent [11] J. Kigami, Harmonic metric and Dirichlet form on the Sierpinski gasket, in K. D. Elworthy and N. Ikeda (eds), Asymptotic problems in probability theory: stochastic models and diffusion on fractals (Sanda/Kyoto, 1990), Pitman research notes in Math., {\bf 283}, Harlow, 201-218, 1993.

\noindent [12] P. Koskela, Y. Zhou, Geometry and Analysis of Dirichlet forms, Adv. Math., 
{\bf 231}, 2755-2801, 2012. 


\noindent [13] S. Kusuoka, Dirichlet forms on fractals and products of random matrices, Publ. Res. Inst. Math. Sci., {\bf 25}, 659-680, 1989.

\noindent [14] R. Ma\~n\'e, Ergodic theory and differentiable dynamics, Berlin, 1983.




\noindent [15] R. Peirone, Convergence of Dirichlet forms on fractals, mimeographed notes.

\noindent [16] W. Perry, M. Pollicott, Zeta functions and the periodic orbit structure of hyperbolic dynamics, Asterisque, {\bf 187-188}, 1990. 

\noindent [17] V. Rokhlin, On the fundamental ideas of measure theory, Transl. Amer. Math. Soc., {\bf 71}, 1952. 

\noindent [18]  Disintegration into conditional measures: Rokhlin's theorem, notes on M. Viana's website. 

\noindent [19] A. Teplyaev, Energy and Laplacian on the Sierpi\`nski gasket, Proceedings of symposia in pure Mathematics, {\bf 72-1}, 131-154, 2004. 

\noindent [20] M. Viana, Stochastic analysis of deterministic systems, mimeographed notes.

\end

By Fubini, there is $A\subset[-1,1]$ with $\L^1(A^c)=0$ such that every $t\in A$ is a point of differentiability of $\g$ for $\pi$-a. e. $\g\in C([-1,1],S)$; calling $B_t$ this full measure set, it suffices to show (4.6) when $t\in A$ and $\g\in B_t$. 

By lemma 4.1 we know that 
$$\hat A=\{
t\in A\st\g_t\in F
\}   $$
is again a full measure set; thus, it suffices to show (4.6) for $t\in\hat A$. 

We saw in section 1 that the restriction of a test plan is another test plan; thus, in (4.6) we can as well suppose that $t=0$; we have to show that, if $\g$ is differentiable at $t=0$ and $\g_t\in F$, then 
$$\dot\g_0=\l v_{\g(0)}    \eqno (4.6)$$
for $\pi$-a. e. $\g\in B_t$. 

We show this. Let us disintegrate $\pi$ with respect to the evaluation map $e_0$, i. e. 
$$\pi=[(e_0)_\sharp\pi]\otimes\mu_x  \eqno (4.6)$$
where $\mu_x$ is a Borel measure on $C([-1,1],S)$ which concentrates on the curves with 
$\g_0=x$. By (4.6) it suffices to show that, for $(e_0)_\sharp\pi$-a. e. $x\in \hat F$ and 
$\mu_x$-a. e. $\g\in C([-1,1],S)$, we have (4.6). Since $\pi$ has bounded deformation, we have that $(e_t)_\sharp\pi<<\kappa$; thus, it suffices to prove (4.6) for $\kappa$-a. e. 
$x\in F$ and $\mu_x$-a. e. $\g\in C([-1,1],S)$. Since this is the statement of lemma 4.1, we are done.